\documentclass[final,11pt,3p]{elsarticle}

\usepackage{graphicx}
\usepackage{amssymb}
\usepackage{amsthm}

\newtheorem{theorem}{Theorem}

\newtheorem{definition}{Definition}
\newtheorem{remark}{Remark}
\newtheorem{corollary}{Corollary}
\newtheorem{example}{Example}

\usepackage{lineno}

\makeatletter
\def\ps@pprintTitle{%
 \let\@oddhead\@empty
 \let\@evenhead\@empty
 \def\@oddfoot{}%
 \let\@evenfoot\@oddfoot}
\makeatother

\journal{ICSM-5}

\begin{document}

\begin{frontmatter}


\title{Extremes of Sums and Maxima with Application to Random Networks}



\author[label1]{Natalia Markovich\corref{cor1}}
\cortext[cor1]{Corresponding author}
\ead{nat.markovich@gmail.com}
\address[label1]{V.A. Trapeznikov Institute of Control Sciences,
Russian Academy of Sciences, Moscow, Russian Federation}

\begin{abstract}
The sums and maxima of non-stationary random length sequences of regularly varying random variables   may have the same tail and extremal indices, Markovich and Rodionov (2020).
 The main constraint  is that there exists a unique series in a scheme of series with the minimum tail index. The result is now revised allowing  a
random bounded number of series to have the minimum tail index.
This new result is applied to random networks.
\end{abstract}


\end{frontmatter}
\pagestyle{empty}

\section{Introduction}\label{S:1}
Random length sequences and distribution tails of their sums and maxima attract the interest of many researchers due to numerous applications including  queues, branching processes and random networks \cite{AsmussenFoss2018}, \cite{JesMik2006}, \cite{Lebedev2015a}, \cite{Lebedev2015b}, \cite{MarkRod2020}, \cite{OlveraCravioto2012}, \cite{RobertSegers}, \cite{TillierWintenberger}. \\
Let  $\{Y_{n,i}: n,i\ge 1\}$ be a doubly-indexed array of nonnegative random variables (r.v.s) in which the "row index" $n$ corresponds to time, and the "column index" $i$  enumerates the series. On the same probability space, the existence of a sequence of non-negative integer-valued r.v.s $\{N_n: n\ge 1\}$ is assumed.
Let $\{Y_{n,i}: n\ge 1\}$ be a strict-sense stationary sequence  with extremal index $\theta_i$ having a regularly varying tail 
\begin{eqnarray}\label{11a} P\{Y_{n,i}>x\}&=&\ell_i(x)x^{-k_i}\end{eqnarray}
with tail index $k_i>0$ and  a slowly varying function $\ell_i(x)$. There are no assumptions on the dependence structure in $i$.
In \cite{MarkRod2020}, the weighted sums and maxima
\begin{eqnarray*}\label{3}
&&Y_{n}^*(z, N_n) =\max(z_1Y_{n,1},...,z_{N_n}Y_{n,N_n}),
~~Y_{n}(z, N_n)=z_1Y_{n,1}+...+z_{N_n}Y_{n,N_n}\nonumber\end{eqnarray*}
for positive constants $z_1, z_2,...$  were considered. A similar result was obtained in \cite{Gold:13} for random sequences of a fixed length $l\ge 1$ and when $\{Y_{n,i}: n\ge 1\}$ has a power-type tail, i.e. $P(Y_{n,i}>x) \sim c^{(i)} x^{-k_i}$ as $x\to\infty$, where $c^{(i)}$ is a real-valued positive constant.
\begin{definition}\label{Def1} A stationary sequence  $\{Y_n\}_{n\ge 1}$ with distribution function $F(x)$ and $M_n=\bigvee_{j=1}^{n}Y_j=\max_j Y_j $ is said to have extremal index $\theta\in[0,1]$ if
for each $0<\tau <\infty$ there is a sequence of real numbers $u_n=u_n(\tau)$ such that
\begin{eqnarray*}\label{1}&&\lim_{n\to\infty}n(1-F(u_n))=\tau \qquad\mbox{and}\end{eqnarray*}
\begin{equation}\label{2}\lim_{n\to\infty}P\{M_n\le u_n\}=e^{-\tau\theta}\end{equation}
hold (\cite{Leadbetter}, p.63).
\end{definition}
I.i.d. r.v.s $\{Y_n\}$ give $\theta=1$. The converse may be incorrect.
An extremal index that is close to zero implies a kind of a strong local dependence. 
\\
Let us recall Theorem 4 derived in \cite{MarkRod2020}. It is assumed that the "column" sequences $\{Y_{n,i}: i\ge 1\}$
have stationary distributions (\ref{11a})
 in $n$ with positive tail indices $\{k_1, k_2,...\}$ and
extremal indices $\{\theta_1, \theta_{2},...\}$ for each fixed $i$, where 
$\{\ell_i(x)\}$ are restricted by the 
condition: for all
$A>1$, $\delta>0$
there exists $x_0(A, \delta)$ such that for all $i\geq 1$
\begin{eqnarray} \ell_i(x)\leq A x^\delta,\ \  x>x_0(A, \delta) \label{uniform}
\end{eqnarray}
holds. 
$N_n$ has a regularly varying distribution
 with tail index $\alpha>0,$ that is
 \begin{eqnarray}\label{15a} &&P(N_n>x) = x^{-\alpha} \tilde{\ell}_n(x).\end{eqnarray}
 There is a minimum tail index $k_1$ and
$k:= \lim_{n\to\infty} \inf_{2\leq i\leq l_n} k_i,$ 
 \begin{eqnarray}\label{27}&& l_n=[n^\chi],
  \qquad
   \end{eqnarray}
  and  $\chi$ satisfies
    \begin{equation}
0<\chi<\chi_0, \qquad\chi_0 = \frac{k-k_1}{k_1(k+1)}.
\label{chi}
\end{equation}
An arbitrary dependence structure  between  $\{Y_{n,i}\}$ and  $\{N_n\}$  is allowed. The tail of  $N_n$ does not dominate the tail of the most heavy-tailed term $Y_{n,1}$. Let $u_n=yn^{1/k_1}\ell_1^{\sharp}(n)$, $y>0$, where $\ell_1^{\sharp}(n)$ is the de Brujin conjugate of $\ell(x)=(\ell_1(x))^{-1/k_1}$, and the positive weights $\{z_i\}$ are assumed to be bounded as in \cite{MarkRod2020}.
 \begin{theorem}\label{T1} \cite{MarkRod2020}
Let the sets of slowly varying functions $\{\tilde{\ell}_n(x)\}_{n\geq 1}$ in (\ref{15a}) and $\{\ell_i(x)\}_{i\geq 1}$ in (\ref{11a}) satisfy the condition (\ref{uniform}).
Suppose that $k_1<k$ and
\begin{eqnarray}\label{4a}P\{N_n>l_n\}&=&o\left(P\{Y_{n,1}>u_n\}\right), ~~n\to\infty\end{eqnarray}
hold, where the sequence $l_n$ satisfies (\ref{27}) and (\ref{chi}).
Then the sequences $Y_{n}^*(z,N_n)$ and $Y_{n}(z,N_n)$
have the same tail index
$k_1$ and the same extremal index $\theta_1$.
\end{theorem}
Theorem \ref{T1} is based on Theorem \ref{T2}  we recall further. Let us denote $Y_n^*(z)=Y_n^*(z,l_n)$ and $Y_n(z)=Y_n(z,l_n)$.
\begin{theorem}\label{T2} \cite{MarkRod2020}
Let $k_1<k$, (\ref{uniform}), (\ref{27})
and (\ref{chi}) hold. Then the sequences $Y_n^*(z)$ and $Y_n(z)$ have the same tail index
$k_1$ and the same extremal index $\theta_1$.
\end{theorem}
Our objectives are twofold. At first, we  revise Theorem \ref{T1} for the case when a random number of "column" series may have the minimum tail index.  At second, we  modify Theorem \ref{T1}  with regard to random networks. Each node pair in a random network is connected with some probability \cite{Barabasi}. 
\\
Let $G_n=(V_n, E_n)$ be a directed graph with a set of vertices $V_n=\{1,...,n\}$, and a set of directed edges $E_n$. Google's PageRank  vector $R = (R_1, . . . , R_n)$
is the unique solution to the following system of linear equations:
\begin{eqnarray}\label{00}
R_i&=&c\sum_{j:(j,i)\in E_n}\frac{R_j}{D_j}+(1-c)q_i,\qquad i=1,...,n,
\end{eqnarray}
where 
the summation is taken over a number of pages $j$ that link to page $i$ 
(in-degree), $D_j$ is the number of outgoing links of page $j$ 
(out-degree), $c\in(0,1)$  is a damping factor, $q=(q_1,q_2,...,q_n)$ is a personalization probability vector or user's preferences such that $q_i\ge 0$ and $\sum_{i=1}^nq_i=1$, and $n$ is the total number of pages. 
The World Wide Web (Web) is a very large interconnected graph where nodes correspond to pages.
The PageRank was designed to rank pages on the Web in such a way that a page is important if many important pages have a
hyperlink to it \cite{Chen2014RankingAO}.
\\
A stochastic approach to analyze (\ref{00}) is the following. The PageRank of a randomly chosen Web page (i.e. a vertex on a Web graph) considered as a root node of a Galton-Watson tree with random in- and out-degrees may be modeled  as a r.v. $R$ which is  the solution
of the fixed-point problem (cf.  \cite{OlveraCravioto2012}, \cite{Vol:10}, \cite{ChLiOl:14}, \cite{Jel:10})
 \begin{equation}\label{6} R=^D\sum_{j=1}^{N}A_{j}R_{j}+Q 
 \end{equation}
 assuming that $\{R_j\}$ are independent identically distributed (i.i.d.) copies of $R$ and $E(Q)<1$ holds.
 $(Q,N,\{A_j\})$ is a real-valued vector.
 $N$ denotes the in-degree. $Q$ is a personalization value of the vertex \cite{Vol:10}.
$=^D$ means equality in distribution.
Under the assumptions (we shall call it Assumptions A)
 that
 $\{R_j\}$ are i.i.d. and independent of $(Q, N, \{A_j\})$ with $\{A_j\}$ independent of $(N,Q)$, and that $N$ and $Q$ are allowed to be dependent,   it is stated in \cite{Vol:10}, \cite{Jel:10} that the stationary distribution of $R$ in (\ref{6}) is regularly varying and its tail index is determined by the most heavy-tailed distributed term in the regularly varying distributed pair $(N, Q)$.
  \\
 This result was generalized by \cite{AsmussenFoss2018}, and the unique solution of (\ref{6}) is proved to be intermediate regularly varying\footnote{The class of intermediate regularly varying distributions such that $\lim_{\alpha\uparrow 1}\limsup_{x\to\infty}\overline{F}(\alpha x)/\overline{F}(x)=1$ includes regularly varying distributions.}, if $Q$ or $N$ has an intermediate regularly varying distribution, or $(Q,N)$ has a two-dimensional regularly varying distribution. The multivariate version of (\ref{6})
 \begin{eqnarray*}&& R(i)=^D\sum_{k=1}^{K}\sum_{m=1}^{N^{(k)}(i)}R_{m}(k)+Q(i),
 \end{eqnarray*}
 where $R_{m}(k) =^D R(k)$ holds,
  is considered with similar assumptions and regularly varying statements by thinking that $N^{(k)}(i)$ is a number of type-$k$ children of a type-$i$ ancestor and considering a multi-type Galton-Watson tree.
\\
A Max-Linear Model  \cite{GisKlu:15}  is obtained by substitution sums in  (\ref{00})
by maxima, i.e.
\begin{eqnarray}\label{7}R(i)&=&\bigvee_{j\to i}A_{j}R(j)\vee Q_i,
~~ i=1,...,n.\end{eqnarray}
 Under Assumptions A 
 the power law tail  $P\{|R|>x\}\sim Hx^{-\alpha}$,\footnote{The symbol $\sim$ means asymptotically equal to or $f(x)\sim g(x)$ $\Leftrightarrow$ $f(x)/g(x)\rightarrow 1$ as $x\rightarrow a$, $x\in M$ where the functions  $f(x)$ and $g(x)$ are defined on some set $M$ and $a$ is a limit point of $M$.} $\alpha>0$, $H>0$ as $x\to\infty$ of the so called 'minimal/endogeneous' solution
  of the  equation
  \begin{equation}\label{6a}R=^D\left(\bigvee_{j=1}^{N}A_{j}R_j\right)\vee Q,\end{equation}
 is derived in \cite{Jel:15}.
 We propose to apply Theorem \ref{T1} to the right-hand sides of  (\ref{6}) and (\ref{6a}) under weaker assumptions than the Assumptions A. Namely, the conditions that PageRanks $\{R_j\}$ of the first generation are i.i.d. and the mutual independence of $\{R_j\}$ and $N$ are omitted.
 By a modified Theorem \ref{T1} (Theorem \ref{T5}) we obtain that the most heavy-tailed terms $\{R_{j}\}$ determine the heaviness of the tail and the extremal index  of both the sum and maximum.
 \\
The paper is organized as follows. The revision of Theorems \ref{T1} and \ref{T2} by Theorems \ref{T3} and \ref{T4} is given in Section \ref{Sec2_1}. The modification of Theorem \ref{T4} by Theorem \ref{T5} to find  the extremal and tail indices of PageRank and the Max-Linear Model as influence measures of nodes in 
random direct graphs is presented in Section \ref{Sec2_2}. 
Conclusions and a discussion are presented in Section \ref{Sec3}. The proofs are given in Section \ref{Sec4}.
\section{Revision of Theorem \ref{T1}}\label{Sec2_1}

We  revise Theorem \ref{T1} by Theorem \ref{T4}  allowing  a random bounded number $d\ge 1$ of series to have a minimum tail index.
To this end, we extend Theorem \ref{T2} by Theorem \ref{T3}. 
We assume in Theorems \ref{T3} and \ref{T4} that $k_i=k_1$, $i\in\{1,..., d\}$,  $1\le d\le l_n-1$,
$k_1<k$, where
\begin{eqnarray}\label{5}k &:=& \lim_{n\to\infty} \inf_{d+1\leq i\leq l_n} k_i,\end{eqnarray}
and that 
(\ref{27}), (\ref{chi}) hold.
\\
We introduce the following independence condition:
\\
(A1) The stationary sequences $\{Y_{n,i}\}_{n\ge 1}$, $i\in\{1,...,d\}$ are mutually independent, and  independent of the sequences
$\{Y_{n,i}\}_{n\ge 1}$, $i\in\{d+1,...,l_n\}$. 
\\
Denote $M_n^{(i)}= \max\{Y_{1,i}, Y_{2,i},...,Y_{n,i}\}, ~i\in\{1,..,l_n\}$.
\begin{theorem}\label{T3} Let
(\ref{uniform}) hold for all $d+1\leq i\leq l_n$.
Then the sequences
$Y_n^*(z,l_n)$ and $Y_n(z,l_n)$ have the same tail index
$k_1$.
\begin{enumerate}
\item If, in addition, (A1) holds, then
 $Y_n^*(z,l_n)$ and $Y_n(z,l_n)$ have  the same extremal index
\begin{eqnarray}\label{5aaab}\theta(z) &= & \sum_{j=1}^d\theta_jz_j^{k_1}/\sum_{j=1}^dz_j^{k_1}.
\end{eqnarray}
\item If, instead, (A1) does not hold and $d=1$, then
their extremal index is equal to $\theta_1$, but it may 
not exist, if $d>1$. Assuming for $d>1$ that  all elements in pairs of the "column" sequences $\{Y_{n,i}\}_{n\ge 1}$, $i\in\{1,...,d\}$ have the same mutual dependence, and \begin{eqnarray}\label{11b}&&\sum_{j=1}^{d-1} P\{z_{j}M_n^{(j)}>u_n,z_{j+1}M_n^{(j+1)}\le u_n,...,z_{d}M_n^{(d)}\le u_n\}=o(P\{z_dM_n^{(d)}\le u_n\})
\end{eqnarray}
 holds as $n\to\infty$, $Y_n^*(z,l_n)$ and $Y_n(z,l_n)$ have  the  extremal index  $\theta_d$.
 \end{enumerate}
\end{theorem}
\begin{remark}Since an enumeration of the first $d$ "column" sequences is not significant, one can rewrite the condition (\ref{11b}) with regard to the $i$th column and obtain the same
 statement for any $\theta_i$, $i\in[1,d]$.  \end{remark}
\begin{corollary}\label{Cor1} Let $Y_{n,i}=Y_{n,1}$, 
$n\ge 1$, $i\in\{1,...,d\}$. 
Then $Y_n^*(z,l_n)$ and $Y_n(z,l_n)$ have the same tail index
$k_1$ and the same extremal index $\theta_1$.
\end{corollary}
\begin{remark}
If there are in total $d+1$ stationary mutually independent "column" sequences having the same tail index $k_1$ and extremal indices $\theta_1,...,\theta_{d+1}$,
then $Y_n^*(z,l_n)$ and $Y_n(z,l_n)$ have the tail index $k_1$
and the extremal index  that is a superposition of $\theta_1,...,\theta_{d+1}$ 
 as derived in \cite{Gold:13}, see  Theorem 2 in \cite{MarkRod2020}.
\end{remark}
\begin{example}\label{Exam1} In case of an arbitrary dependence among elements of the $d$
"columns" series with the minimum tail index,
the extremal index of maxima $Y_n^*(z,l_n)$ and sums   $Y_n(z,l_n)$ may not exist.
Suppose elements of odd "row" sequences coincide and elements of the even rows are i.i.d..
 Then sums and maxima over rows are differently distributed and
the sequences $Y_n^*(z,l_n)$ and $Y_n(z,l_n)$ are nonstationary.
\end{example}
\begin{example}\label{Exam2} The assumption (\ref{11b}) is valid for the following $d$ "column" sequences. Let elements of each  "column" sequence be sums 
 of corresponding elements of all previous "columns", i.e. $Y_{n,i}=\sum_{j=1}^{i-1}Y_{n,j}$
 and $z_1\le z_2\le ...\le z_d$. Each of $d$ "column" sequences has the tail index $k_1$ that follows from the first statement of Theorem \ref{T3}. Since $M_n^{(1)}=M_n^{(2)}< M_n^{(3)}<...<M_n^{(d)}$ holds, then  $\sum_{j=1}^{d-1} P\{z_{j}M_n^{(j)}>u_n,z_{j+1}M_n^{(j+1)}\le u_n,...,z_{d}M_n^{(d)}\le u_n\}=0$ holds.
\end{example}
Now we reformulate Theorem \ref{T1}.
\begin{theorem}\label{T4}
Let the sets of slowly varying functions $\{\tilde{\ell}_n(x)\}_{n\geq 1}$ in (\ref{15a}) and $\{\ell_i(x)\}_{d+1\le i\le l_n}$ in (\ref{11a}) satisfy the condition (\ref{uniform}).
Assume that
$d$ is a bounded discrete r.v.  such that
$d<d_n = \min(C, l_n)$, 
$C>1$ holds, and
  $d$ and $\{Y_{n,i}\}$ are mutually independent. Let (\ref{4a}) hold.
Then  $Y_{n}^*(z,N_n)$ and $Y_{n}(z,N_n)$
have the same tail index
$k_1$, but their extremal indices do not exist.
\end{theorem}
\begin{corollary}\label{Cor2} The results of Theorems \ref{T3} and \ref{T4}
remain true if the tail indices $\{k_{n,i}\}$ of elements
in the "columns" $\{Y_{n,i}: n\ge 1\}$  are different,
apart of those
columns with the minimum tail
index.
The elements of the columns with non-minimum tail indices 
may be partly light-tailed distributed. 
\end{corollary}
The corollary follows from the proof of Theorem 3 and Lemma 1 in \cite{MarkRod2020}.

\section{Application to Random Networks}\label{Sec2_2}
\begin{figure}[tbp]
 \begin{minipage}[t]{\textwidth}
\centering
\includegraphics[width=0.45\textwidth]{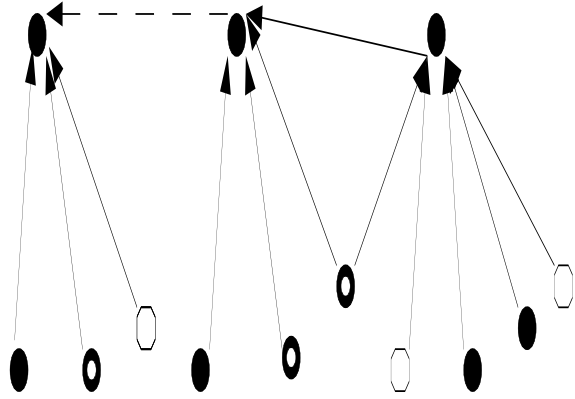}
  \caption{The (disconnected) graph with nodes from different communities (a dashed arrow may be absent). Nodes from dominating communities with the most heavy-tailed distributed influence indices having the minimum tail index
  $k_1$ are  marked by filled circles, and nodes from other communities by open circles and filled rings. 
}
\label{fig:4}
\end{minipage}
\end{figure}

Each node in a random network
may be considered as a root of some 
directed graph of its 
followers (i.e. nodes with incoming links to the root node) which may contain cycles. Thus, generations of followers may be overlapping. The $r$th generation implies the set of nodes  in the tree at distance $r$ from the root node \cite{Chen2014RankingAO}. There may exist links between nodes within generations.
Influence characteristics of nodes such as PageRanks or the Max-Linear Models calculated by these generations
may be dependent.
\\
Generations  may be non-stationary distributed since their nodes 
may belong to 
communities with different distributions  and the graph is not necessarily fully connected, Fig.\ref{fig:4}. 
A community
structure is  'the organization of vertices in clusters, with many edges joining
vertices of the same cluster and comparatively few edges joining vertices of different clusters' \cite{Fortunato}. There are methods to extract the community structure of large networks, see, for instance, \cite{Blondel_2008}, \cite{Newman}.
\\
Theorem \ref{T4} is further reformulated in the context of  PageRank and the Max-Linear Model. Denote  an
in-degree of the node $i$ as $N_i$. $N_i$ determines the random size of the generation   of a one-link neighborhood from the node $i$.
Let us denote 
 $A_jR(j)=cR(j)/D_j$,
$j\in\{1,..., N_i\}$ in (\ref{00}) and (\ref{7}) as $z_jY_{i,j}$ with $z_j=c$.
Then one can rewrite the right-hand sides of (\ref{6}) and (\ref{6a}) in notations of Theorems \ref{T1} and \ref{T4}
as
\begin{eqnarray}\label{9a}Y_i(c, N_i)&=&c\sum_{j=1}^{N_i}Y_{i,j}+Q_i,
~~Y^*_i(c, N_i)=c\bigvee_{j=1}^{N_i}Y_{i,j}\vee Q_i,~~ i\in\{1,...,n\},
\end{eqnarray}
where $\{Y_{i,j}\}$ are not i.i.d. r.v.s. $\{Y_{i,j}\}$ and $Q_i$ as well as $\{Y_{i,j}\}$ and $N_i$ are not necessarily mutually independent. $N_i$ satisfies (\ref{4a}).
In the context of PageRank $Q_i=(1-c)q_i$ holds.
\\
Let us consider the matrix (\ref{4}) of the
array $\{Y_{n,i}: n,i\ge 1\}$ corresponding to (\ref{9a}),
 and the tail and extremal indices of its columns are shown in matrix (\ref{4b}):
\begin{equation}\label{4}\left(
    \begin{array}{ccccccc}
      cY_{1,1} & cY_{1,2} & cY_{1,N_1}&...  & 0 & 0 & Q_1 \\
      cY_{2,1} & 0 & cY_{2,3} & ~~...  & cY_{n,N_2} & 0 & Q_2 \\
      ...& ... & ... & ... & ... & ... \\
     0 & cY_{n,2} & cY_{n,3} & ~~... & 0 & cY_{n,N_n} & Q_n \\
    \end{array}
  \right),
\end{equation}
\begin{equation}\label{4b}
\left(
    \begin{array}{ccccccc}
     k_1~~~ & k_2 ~~~& k_{3}~~~ & ... ~~ & ... &~~~~ k_{N_n}  & ~\beta \\
     \theta_1~~~ & \theta_2 ~~~& \theta_{3}~~~ & ... ~~& ...  &  ~~~~\theta_{N_n} & 1 \\
          \end{array}
  \right).
\end{equation}
The sequence $\{N_n\}$ is equal to the number of non-zero elements in the rows and it is not necessarily increasing. 
Row elements in (\ref{4}) (excluding $\{Q_i\}$) correspond to the first generations of followers 
of  nodes that can be arbitrary enumerated as $1,2,...,n$.
The column $i$ contains a community 
and  has the extremal index $\theta_i$. Zeroes in the rows imply that the
first generations do not contain followers from the corresponding communities.
"Column" series 
may be mutually dependent and have tails (\ref{11a}) 
with 
tail indices $\{k_i\}$. 
The latter columns  may (partly) coincide.
  In random networks the column dependence reflects the fact that 
  communities may be overlapping.
 \\
 The last column in (\ref{4}) contains i.i.d. r.v.s of user preferences $\{Q_i\}$.
 Theorem \ref{T4} is  valid if $\{Q_i\}$ are heavy- or light-tailed distributed.
  In- and out-degrees of nodes are  available statistics in practice, but the user preference has to be simulated \cite{Vol:10}. The preference is often expected to be uniformly distributed, i.e. an arbitrary page can be selected by a user uniformly among the nodes of the network. In \cite{Vol:10} Pareto distributed 
  $\{Q_i\}$ were studied by real data.
 The tail of $\{N_i\}$ is assumed here to be lighter than the regularly varying tails of $\{Y_{i,j}\}$ and $\{Q_i\}$. We assume that there
  exist  columns
 with a minimum tail index $k_1$, $k_1<k$.
 \\
Let the graph 
be partitioned into 
communities of nodes whose 
PageRanks are stationary regularly varying distributed with tail indices $\{k_1, k_2,...\}$ and extremal indices $\{\theta_1, \theta_2,...\}$. 
We assume that at least one neighbor of a root-node $i$ belongs to  one of the $d$ most heavy-tailed (dominating) communities (i.e. there exists a type-$k_1$ child of the ancestor $i$). The number of such neighbors is random and they may belong to different dominating communities. 
 It means that all row sequences $\{Y_{n,i}:n\ge 1\}$ in (\ref{4}) must have at least one 
element with the smallest tail index.
\begin{theorem}\label{T5}
PageRanks and the Max-Linear Models of $n$ root-nodes  
have the same minimum tail index as one of the $d$ most heavy-tailed distributed communities
within a random 
graph associated with the roots 
as $n\to\infty$.
Their extremal index is equal to the extremal index of the most heavy-tailed distributed community if the latter is unique, and it is calculated by
(\ref{5aaab}) using extremal indices of the $d$ most heavy-tailed communities
if the latter are mutually independent and independent of the rest of communities,
or it may 
not exist 
if the latter independence condition is not valid.
\end{theorem}
Theorem \ref{T5} determines the extremal index of the node in the directed graph (except for the graph leaves).
\begin{remark} There are three important practical assumptions in Theorem \ref{T5}.
\\
(1) By Corollary \ref{Cor2} the node influence indices within non-dominating communities are allowed to be non-stationary distributed.
This property allows us to use non-stationary distributed 
communities if they are not tail dominant.
\\
(2) All elements of the dominating communities have to be regularly varying distributed. 
To our best knowledge, PageRanks are proved  to be power law distributed
for branching trees in \cite{Vol:10},  the directed configuration model in \cite{LeeOlveraCravioto}
 and directed generalized random graphs in \cite{ChLiOl:14} under Assumptions A.
 The existence of an asymptotic PageRank distribution for directed graphs
 is proved assuming a local weak convergence of a sequence of directed graphs to a limiting graph \cite{Garavaglia}.
 The conditions under which the PageRank of the root in the limiting graph shows a power-law tail was not found in  \cite{Garavaglia}.
 \\
(3) The dependence between communities impacts the extremal index of PageRanks and the Max-Linear Models of the root-nodes.  
\end{remark}
Our approach can be used to simulate an enlargement of graphs with given values of the tail and extremal indices of the nodes.
To this end, one can simulate communities of nodes  with given tail and extremal indices connected to a sequence of  root-nodes which do not belong to the communities.
If there is a unique community with the minimum tail index among them, then the sequence of the root nodes inherits its tail and extremal indices.
Selecting $d$ mutually independent communities with a minimum tail index $k_1$ and extremal indices $\theta_1,...,\theta_d$,
 PageRank and the Max-Linear Model sequences of the associated root nodes will have the same tail index $k_1$ and  
the extremal index $\theta(z)$ calculated by (\ref{5aaab}) using  $\theta_1,...,\theta_d$
and taking weights $z_j=c$, $j\in\{1,...,d\}$. Adding the sequence of roots as a new community does not change the dependence  structure and the heaviness of tail of PageRank and the Max-Linear Model of the next sequence of roots, namely, it will be $(k_1,\theta(z))$.
\\
Since the reciprocal of the extremal index approximates the mean cluster size \cite{Leadbetter}, the extremal index can  also be used as the node influence. Namely, a node with the extremal index close to zero may be considered as
an influential one since there are generations (the clusters) of followers in its coupled 
graph with a large number of nodes whose influence measures exceed a sufficiently high threshold $u$. One can consider the
generations as random length blocks that can be overlapping.
\\
Theorem \ref{T5} is valid both for  directed and undirected graphs. The calculation of PageRank  by undirected graphs is considered in \cite{Avrachenkov}, for instance.
\section{Conclusions and discussion}\label{Sec3}
  This paper makes a fundamental step forward by extending the analysis of
PageRank and the Max-Linear Model to its extremes and to graphs that are not necessarily trees.
Our approach to find the asymptotic distributions and local cluster properties  of PageRank
and the Max-Linear Model is based on results derived for the sums and maxima of non-stationary random lengths sequences in \cite{MarkRod2020}.
  The main constraint of these results is that a unique series with the minimum tail index in a scheme of series with
 regularly varying tails is assumed. This assumption may not be realistic in practice since a node in a random network may have  a random number of
 followers (not necessarily  a unique follower) belonging to the most heavy-tailed communities. Another important assumption is that
 the mutually dependent series are allowed to be regularly varying. Regarding PageRank and the Max-Linear Model 
 the regularly varying assumption of the latter characteristics is crucial, since the power-law tails are proved to our best knowledge
 only for mutually independent series.
  \\
Our results are twofold. At first, we extend  Theorems 3 and 4 in  \cite{MarkRod2020} by Theorems \ref{T3} and \ref{T4}. 
 Secondly, we interpret these results to PageRank and the Max-Linear Model.
  Theorem \ref{T5} is the analogue of Theorem \ref{T4} for random
 networks.
 \\
 To our best knowledge, the extremal index of PageRank and the Max-Linear Model is considered at first time. The extremal index may
 serve  as a new influence measure of nodes.
\section{Proofs}\label{Sec4}
\subsection{Proof of Theorem \ref{T3}}
The proof is similar to that of Theorem 3 \cite{MarkRod2020}.  We just indicate the modifications.  The numeration of Theorem 3 \cite{MarkRod2020} is preserved throughout the proof. Let us take the same sequence of thresholds $u_n = y n^{1/k_1} \ell^{\sharp}_1(n)$, $y>0$ as in \cite{MarkRod2020}.
\paragraph{\textbf{Tail index}}
At first, we show that the tail index of $Y_n(z)$ and $Y_n^*(z)$ is the same. The right-hand side of (12) in \cite{MarkRod2020} can be rewritten as
\begin{eqnarray}\label{14}P\{Y_n(z)>u_n\}&\le & P\{\sum_{i=1}^{d}z_iY_{n,i}>u_n(1-\varepsilon)\}+\sum_{i=d+1}^{l_n}P\{z_iY_{n,i}>u_n\varepsilon_i\},
\end{eqnarray}
where  $\sum_{i=1}^{l_n}\varepsilon_i=1$ holds and $\{\varepsilon_i\}$  is a sequence of positive elements. Let us denote $\varepsilon=\sum_{i=d+1}^{l_n}\varepsilon_i$.
One may take $\varepsilon_i$, $i\in\{d+1,...,l_n\}$, in such a way
to satisfy $\varepsilon_i\to 0$ and $\varepsilon\to 0$ as $n\to\infty$. Moreover, choosing $\{\varepsilon_i=1/l_n^{\eta+1}\}$, $\eta>0$ as in Lemma 1 in  \cite{MarkRod2020} one can derive (13) in \cite{MarkRod2020} substituting $2$ by $d+1$, namely, the following
\begin{eqnarray}\label{5a}&&\sum_{i=d+1}^{l_n}P\{z_iY_{n,i}>u_n\varepsilon_i\}=o(1/n)\end{eqnarray}
as $n\to\infty$. To prove the latter, we need to assume (\ref{uniform}) for $d+1\leq i\leq l_n$. For simplicity, let us consider the case $d=2$.
Then by (7) in \cite{MarkRod2020}
\begin{eqnarray}P\{z_1 Y_{n,1} > u_n\}  = (z_1/y)^{k_1} n^{-1} (1+o(1)), ~~n\to\infty\label{quantile}\end{eqnarray}
and (13) in \cite{MarkRod2020} we obtain
\begin{eqnarray*}&&P\{z_1Y_{n,1}+z_2Y_{n,2}>u_n(1-\varepsilon)\}
\\
&\le &P\{z_1Y_{n,1}>u_n(1-\varepsilon)\varepsilon_1^*\}+P\{z_2Y_{n,2}>u_n(1-\varepsilon)\varepsilon_2^*\}
\\
&=& \frac{n^{-1}}{(y(1-\varepsilon))^{k_1}}\left[\left(\frac{z_1}{\varepsilon_1^*}\right)^{k_1}+\left(\frac{z_2}{\varepsilon_2^*}\right)^{k_1}\right](1+o(1)),
\end{eqnarray*}
where $\varepsilon_1^*+\varepsilon_2^*=1$ holds. Let us denote
\begin{eqnarray*}(z^*)^{k_1}&=& \left(\frac{z_1}{\varepsilon_1^*}\right)^{k_1}+\left(\frac{z_2}{\varepsilon_2^*}\right)^{k_1}. 
\end{eqnarray*}
By (14) in \cite{MarkRod2020} it follows
\begin{eqnarray}\label{8a}P\{Y_n(z)>u_n\}&\le & \left(\frac{z^*}{y(1-\varepsilon)}\right)^{k_1}n^{-1}(1+o(1)) + o(1/n).
\end{eqnarray}
From another side, we have
\begin{eqnarray}\label{8b}P\{Y_n(z)>u_n\}&\ge & 
P\{z_1Y_{n,1}>u_n\}=
\left(\frac{z_1}{y}\right)^{k_1}n^{-1}(1+o(1))
\end{eqnarray}
 and the left-hand side of (9) in \cite{MarkRod2020} does not change. The heaviness of tail of both  $Y_n(z)$ and $Y_n^*(z)$ coincides.
The  proof is the same for $d>2$.
\\
\paragraph{\textbf{Extremal index for $d$ independent "column" sequences}}
We  assume  the condition (A1).
Denoting
\begin{eqnarray*}
&&M_n^*(z)= \max\{Y_1^*(z),...,Y_n^*(z)\}=\max\{z_1M_n^{(1)},...,z_{l_n}M_n^{(l_n)}\}, ~n\ge 1,
\end{eqnarray*}
 where $Y^*_{n}(z)=Y^*_{n}(z,l_n)$ 
 as in (18) \cite{MarkRod2020}, we have due to independence
\begin{eqnarray*}
&&P\{M_n^*(z)\le u_n\} = \prod_{j=1}^{d}P\{z_j M_n^{(j)}\le u_n\}P\{z_{d+1}M_n^{(d+1)}\le u_n,...,z_{l_n}M_n^{(l_n)}\le u_n\}.
\end{eqnarray*}
Since the extremal index of the sequence $\{Y_{n,j}\}_{n \geq 1}$, $1\le j\le d$  
is assumed to be equal to $\theta_j,$ by (\ref{2}) and (\ref{quantile})
 we get
\begin{eqnarray*}\lim_{n\to\infty}\prod_{j=1}^{d}P\{z_jM_n^{(j)}\le u_n\}&=& \exp(-\sum_{j=1}^{d}\theta_j(z_j/y)^{k_1}). \end{eqnarray*}
By (19) and (21) in \cite{MarkRod2020} and by  (\ref{2}), we obtain
\begin{eqnarray*}
&& P\{z_{d+1}M_n^{(d+1)}\le u_n,...,z_{l_n}M_n^{(l_n)}\le u_n\}= P\{z_{d+1} M_n^{(d+1)} \leq u_n\} (1+o(1))\to 1,
\end{eqnarray*}
since
\begin{eqnarray*}
&& nP\{z_{d+1}Y_{n,d+1}>u_n\}\sim n\cdot n^{-k_{d+1}/k_1}\ell(n)\to 0
\end{eqnarray*}
as $n\to\infty$ due to $k_{d+1}>k_1$.
We obtain
\begin{eqnarray}\label{10}
&& \lim_{n\to\infty}nP\{Y_n^*(z)> u_n\}=\sum_{j=1}^d\left(\frac{z_{j}}{y}\right)^{k_1}
\end{eqnarray}
since it holds
\begin{eqnarray*}
&&P\{Y_n^*(z)> u_n\}=P\{\max(z_1Y_{n,1},...,z_{d}Y_{n,d},z_{d+1}Y_{n,d+1},...,z_{l_n}Y_{n,l_n})> u_n\}
\\
&\sim & \!\!\!\sum_{i=1}^dP\{z_i Y_{n,i}> u_n\}
+P\{\max(z_{1}Y_{n,1},...,z_{d}Y_{n,d})\le u_n\}P\{\max(z_{d+1}Y_{n,d+1},...,z_{l_n}Y_{n,l_n})> u_n\}
\\
&=& \sum_{j=1}^d\left(\frac{z_{j}}{y}\right)^{k_1}\cdot n^{-1}(1+o(1))+o(1/n) 
\end{eqnarray*}
due to (9), (12), (13) and Lemma 1 all in \cite{MarkRod2020}.
Then from \begin{eqnarray}\label{16}
&&P\{M_n^*(z)\le u_n\} =\exp(-\sum_{j=1}^{d}\theta_j(z_j/y)^{k_1}) (1+o(1))
\end{eqnarray}
 we obtain that the extremal index of $Y_n^*(z)$ is equal to (\ref{5aaab}). In the same way as in the proof of Theorem 3 in \cite{MarkRod2020} one can show that $Y_n(z)$ has the same extremal index.
 \paragraph{\textbf{Extremal index for $d$ arbitrary 
 dependent "column" sequences}}
We show that the extremal indices of the sequences $Y_n^*(z)$ and $Y_n(z)$ may not exist if the condition (A1) is not valid. If $d=1$, i.e. there is a unique "column" series with the minimum tail index, we are in the conditions of Theorem \ref{T2}. Let us consider the case $d>1$.
\\
Note, that $Y_n^*(z)$ and $Y_n(z)$ may be non-stationary distributed if the pair-wise dependence between elements of  "column" sequences is different as in Example \ref{Exam1}.
Assume further that  $Y_n^*(z)$ and $Y_n(z)$ are stationary distributed.
Let us rewrite (19)  in \cite{MarkRod2020} as
\begin{eqnarray}\label{12}
&&P\{M_{1,d}\le u_n\}  -  \sum_{i=d+1}^{l_n}P\{z_i M_n^{(i)}> u_n\}\nonumber
\\
&\leq & P\{z_1M_n^{(1)}\le u_n,...,z_{l_n}M_n^{(l_n)}\le u_n\} = P\{M_n^*(z)\le u_n\} \nonumber
\\
&\leq & P\{M_{1,d}\le u_n\}\leq P\{z_dM_n^{(d)}\le u_n\},
\end{eqnarray}
where 
$M_{1,d} = \max\{z_1M_n^{(1)}, z_2M_n^{(2)},...,z_d M_n^{(d)}\}$.
By (21) in \cite{MarkRod2020} and assuming (\ref{uniform}), we have
\begin{eqnarray}\label{12a}\sum_{i=d+1}^{l_n}P\{z_i M_n^{(i)}> u_n\}&=& o(1), ~~n\to\infty.
\end{eqnarray}
We get
\begin{eqnarray*}
&&P\{M_{1,d}\le u_n\}=1-P\{M_{1,d}> u_n\}
\\
&= & 1-P\{z_dM_n^{(d)}> u_n\}-\sum_{j=1}^{d-1} P\{z_{j}M_n^{(j)}>u_n,z_{j+1}M_n^{(j+1)}\le u_n,...,z_{d}M_n^{(d)}\le u_n\}.\nonumber
\end{eqnarray*}
Assuming that (\ref{11b})
 holds, 
 the expression
\begin{eqnarray*}\label{22ab}
P\{M_n^*(z) \leq u_n\} & = & P\{z_d M_n^{(d)} \leq u_n\}(1+o(1))
=\exp(-\theta_d (z_d/y)^{k_1})(1+o(1))
\end{eqnarray*}
that is required to obtain the extremal index $\theta_d$ for $Y_n^*(z)$ follows, otherwise not. The same  holds for $Y_n(z)$ since
\begin{eqnarray}\label{22abc}P\{M_n^*(z)\le u_n\} = P\{M_n(z)\le u_n\} (1+o(1))\end{eqnarray}
can be derived similarly to (23) in \cite{MarkRod2020} due to the common tail index of $Y_n^*(z)$ and $Y_n(z)$.
\subsection{Proof of Corollary \ref{Cor1}}
The existence of a common tail index $k_1$ for $Y_n^*(z,l_n)$ and $Y_n(z,l_n)$ can be shown the same way as in the proof of Theorem \ref{T3}.
Since $Y_n^*(z,d)=\bigvee_{i=1}^d z_i Y_{n,1}=z^{**}Y_{n,1}$ and $Y_n(z,d)=\sum_{i=1}^d z_i Y_{n,1}=z^{*}Y_{n,1}$ hold and by (\ref{14}), (\ref{5a}), (\ref{quantile}) we get
\begin{eqnarray*}\label{17}&&P\{Y_n(z,l_n) > u_n\}  = 
 \left(z^{*}/y\right)^{k_1} n^{-1} (1+o(1)),\end{eqnarray*}
as $n\to\infty$. The same is valid for $Y_n^*(z,l_n)$ since $z^{**}<z^*$. Due to $M_{1,d}=M_n^{(1)}z^{**}$ and by (\ref{12}), (\ref{12a}) we obtain
\begin{eqnarray*}&&P\{M_n^*(z)\le u_n\}=P\{M_n^{(1)}z^{**}\le u_n\}(1+o(1))
=\exp(-\theta_1(z^{**}z_1/y)^{k_1})(1+o(1))\end{eqnarray*} due to (\ref{quantile}). The same is valid for $M_n(z)$ as in Theorem 3 in \cite{MarkRod2020}.
Thus, the statement follows. 

\subsection{Proof of Theorem \ref{T4}}
Denote $S_{n,d}=\sum_{i=1}^{d}z_iY_{n,i}$ and $S_{n,l_n-d}=\sum_{i=d+1}^{l_n}z_iY_{n,i}$. Since $d$ is bounded
by $d_n=\min(C,l_n)$
 we get
\begin{eqnarray*}\label{13aa}&&P\{Y_n(z)>u_n\}=P\{S_{n,d}+S_{n,l_n-d}>u_n\}\nonumber
\\
&=& 
P\{S_{n,d}+S_{n,l_n-d}>u_n, 1\le d\le \lfloor d_n-1\rfloor\}+P\{S_{n,d}+S_{n,l_n-d}>u_n, d> d_n-1\}\nonumber
\\
&= & \sum_{m=1}^{\lfloor d_n-1\rfloor}P\{S_{n,m}+S_{n,l_n-m}>u_n\}P\{d=m\}+o(1),~~n\to\infty. 
\end{eqnarray*}
In the same way it follows
\begin{eqnarray*}&&P\{Y_n(z)>u_n\}\ge P\{S_{n,d}>u_n, 1\le d\le \lfloor d_n-1\rfloor\} 
= \sum_{m=1}^{\lfloor d_n-1\rfloor}P\{S_{n,m}>u_n\}P\{d=m\}. 
\end{eqnarray*}
Similarly,  one can obtain the lower bound of $P\{Y_n^*(z)>u_n\}$ replacing sums by maxima.
Applying the same steps of the proof to $P\{S_{n,m}+S_{n,l_n-m}>u_n\}$ and $P\{S_{n,m}>u_n\}$  as in  Theorem \ref{T3} we obtain that the sequences
$Y_n^*(z,l_n)$ and $Y_n(z,l_n)$ have the same tail index
$k_1$.  
\\
We investigate the extremal indices of 
$Y_n^*(z,l_n)$ and $Y_n(z,l_n)$ 
in case of the independence condition (A1).
If $d_n=C$ holds, then by (\ref{10}) and (\ref{16}) it follows
\begin{eqnarray*}
&&P\{M_n^*(z)\le u_n\}= \sum_{m=1}^{\lfloor C-1\rfloor}\exp\left(-\sum_{j=1}^m\theta_j(z_j/y)^{k_1}\right) P\{d=m\}(1+o(1))=A(y),
\end{eqnarray*}
\begin{eqnarray*}
&& \lim_{n\to\infty}nP\{Y_n^*(z)> u_n\}=\sum_{m=1}^{\lfloor C-1\rfloor}\sum_{j=1}^m(z_j/y)^{k_1} P\{d=m\}=\tau(y).
\end{eqnarray*}
The expression
\begin{eqnarray*}\label{5aaa}
\!\!\!\theta(z) &=& -\ln\left(A(y)\right)/\tau(y)
\end{eqnarray*}
cannot be considered as an extremal index of $Y_n^*(z)$ due to the presence of the arbitrary constant $y>0$ which is included in $u_n$.
The same result is valid for $Y_n(z)$ due to (\ref{22abc}). 
 The same proof follows for $d_n=l_n<C$
because of the majorant converging series.
\\
If (A1) does not hold, then the extremal indices of $Y_n^*(z)$ and $Y_n(z)$ do 
not exist by similar reasons. 
In case that (\ref{22ab}) is fulfilled, one can indicate only bounds for $P\{Y_n>u_n\}$ and $P\{Y^*_n>u_n\}$ by (\ref{8a}) and (\ref{8b}).
\subsection{Proof of Theorem \ref{T5}}
The statement is based on the proof of 
Theorems \ref{T3} and \ref{T4}. 
According to (\ref{9a}) 
sums $\{Y_i(c,N_i)\}$
and maxima $\{Y_i^*(c,N_i)\}$, $i\in\{1,...,n\}$
determine
PageRanks and the Max-Linear Models of $n$ nodes in the graph, where $Y_i(c,N_i)$ and $Y_i^*(c,N_i)$ are built by PageRanks of neighbors of node $i$ considered as the root of the coupled tree.  
These neighbors belong to some communities.
In the same way as in the proof of Theorem \ref{T4}, one can get that $Y_n(c,N_n)$ and $Y^*_n(c,N_n)$ have the minimum tail index as 
the most heavy-tailed distributed communities, let's say $k_1$, as $n\to\infty$.
Following the proof of 
Theorem 4 in \cite{MarkRod2020}  and considering the maxima of sequences 
$\{Y_i(c,N_i)\}$ and $\{Y^*_i(c,N_i)\}$, $i\in\{1,...,n\}$
one can confirm that these maxima have the same asymptotic distribution as $n\to\infty$ and thus, PageRanks and the Max-Linear Models of $n$ root-nodes 
have the 
extremal index $\theta_1$ of  the most heavy-tailed distributed community in case the latter is unique. If  the number of the  most heavy-tailed distributed communities related to $n$ root-nodes is fixed  and the PageRanks of the neighbors from these communities are mutually independent, then according to Theorem \ref{T3} 
the extremal index of the PageRanks and the Max-Linear Models of the roots is calculated by (\ref{5aaab}). Otherwise, the extremal index of the characteristics of the roots may 
not exist.

\section*{Acknowledgments}
The author was partly supported by Russian Foundation for Basic Research (grant 19-01-00090).












\end{document}